\documentclass[12pt]{amsart} 
\usepackage[mathscr]{eucal}
\usepackage{amsmath,amsfonts}
\parskip=\smallskipamount

\newtheorem{theorem}{Theorem}[section]

\newtheorem{proposition}[theorem]{Proposition}
\newtheorem{corollary}[theorem]{Corollary}
\newtheorem{lemma}[theorem]{Lemma}

\newcommand{\CC}{{\mathbb C}}
\newcommand{\NN}{{\mathbb N}}

\newcommand{\FF}{{\mathbb F}}

\makeatletter
\@addtoreset{equation}{section}
\makeatother

\newcommand{\cA}{{\mathcal A}}
\newcommand{\cB}{{\mathcal B}}

\newcommand{\cD}{{\mathcal D}}

\newcommand{\cF}{{\mathcal F}}
\newcommand{\cG}{{\mathcal G}}
\newcommand{\cH}{{\mathcal H}}
\newcommand{\cK}{{\mathcal K}}
\newcommand{\cL}{{\mathcal L}}

\newcommand{\cN}{{\mathcal N}}
\newcommand{\cO}{{\mathcal O}}
\newcommand{\cP}{{\mathcal P}}
\newcommand{\cR}{{\mathcal R}}

\newcommand{\cU}{{\mathcal U}}
\newcommand{\fI}{{\mathfrak I}}

\newcommand{\lkh}{\cL(\cK,\cH)}
\newcommand{\lhk}{\cL(\cH,\cK)}

\newcommand{\lh}{\cL(\cH)}
\newcommand{\lk}{\cL(\cK)}
\newcommand{\iac}{\mathrm i}

\newcommand{\ra}{\rightarrow}

\newdimen\expt
\def\boxit#1{\setbox0\hbox{$\displaystyle{#1}$}
      \hbox{\lower.4\expt
 \hbox{\lower3\expt\hbox{\lower\dp0
      \hbox{\vbox{\hrule height.4\expt
 \hbox{\vrule width.4\expt\hskip3\expt
      \vbox{\vskip3\expt\box0\vskip2\expt}%
 \hskip3\expt\vrule width.4\expt}\hrule height.4\expt}}}}}}

\begin{document}

\bigskip

\title 
{On L.~Schwartz's Boundedness Condition for Kernels} 

\author{T.~Constantinescu} \author{A.~Gheondea} 

\address{Department of Mathematics,
  University of Texas at Dallas,
  Box 830688 Richardson, TX 75083-0688, U. S. A.}
\email{tiberiu@utdallas.edu}

\address{Institutul de Matematic\u{a}
al Academiei Rom\^{a}ne, C.P.~1-764, 70700 Bucure\c{s}ti, Rom\^{a}nia}
\email{\tt gheondea@imar.ro}

\begin{abstract} In previous works we analysed conditions for linearization 
of hermitian kernels. The conditions on the kernel 
turned out to be of a type considered previously by L.~Schwartz in the 
related matter of characterizing the real linear space generated by positive
definite kernels. The aim of the present note is to 
find more concrete expressions of the Schwartz type conditions: in the
Hamburger moment problem for Hankel type kernels on the free semigroup,
in dilation theory (Stinespring type dilations and Haagerup
decomposability), as well as in multi-variable holomorphy. 
Among other things, we prove that any hermitian holomorphic kernel has 
a holomorphic linearization, and hence that hermitian holomorphic 
kernels automatically satisfy L.~Schwartz's boundedness condition.
\end{abstract}

\maketitle

\section{Introduction}
We analysed in \cite{CG1} conditions under which 
the linearization functor produces a Kre\u{\i}n space
from a hermitian kernel, in the spirit of Kolmogorov type decompositions, 
and subsequently in \cite{CG2} we generalized this construction to kernels
invariant under the action of a semigroup with involution.
We also related these 
constructions with the GNS representations of $*$-algebras, an 
issue of some recent interest in quantum field theory with 
indefinite metric (\cite{AGW}, \cite{MPS}, \cite{Ho}, \cite{St}). 
The conditions on the kernel turned out to be of a type considered previously 
by L.~Schwartz in \cite{Sc}
in the related matter of characterizing those hermitian kernels 
that are in the real linear space generated by positive definite kernels. 
These boundedness conditions 
are rather difficult to be verified, see \cite{AGW},\cite{Ho}, 
and their nature is
quite obscure, see \cite{St}.

The aim of the present note is to find more concrete expressions of 
the L.~Schwartz type conditions. We first note that the invariant
Kolmogorov decomposition has a countarpart in the representation theory of
semigroups with involution on reproducing kernel Kre\u\i n spaces. 
Then, it is explained that the type of invariance that is considered 
in \cite{CG2} can be viewed as a Hankel type condition and we apply 
this to a Hamburger moment problem for the free semigroup on $N$ generators.
This is used in order to show how the Schwartz condition is 
somewhat simplified when it is written for generators of a $*$-algebra.

In Section~4 we first show that the Stinespring
dilation of hermitian linear maps fits into the general scheme of
invariant Kolmogorov decompositions and we make explicit the
connection with completely bounded maps and Wittstock's Theorem
\cite{Wi}. This opens the possibility of defining a class of non-hermitian 
decomposable kernels, that may successfully replace the missing class of 
completely bounded kernels, by using a generalization of Haagerup's
decomposable linear mappings on $C^*$-algebras, cf.\ \cite{Ha}. 
An analog of Paulsen's Dilation Theorem for decomposable kernels 
is obtained in Theorem~\ref{decomp}. 

In Section~5 we show that holomorphic kernels in more than one variable have
Kolmogorov decompositions and hence, that hermitian holomorphic 
kernels automatically satisfy L.~Schwartz's boundedness condition.
In view of the transcription between Kolmogorov decompositions and reproducing
kernel spaces, e.g.\ see Theorem~\ref{repro}, this result is an extension of 
the result of D.~Alpay in \cite{Al} proved for one variable holomorphic 
hermitian kernels.

\section{Preliminaries}  

We briefly review the structure of the Kolmogorov decomposition of invariant 
hermitian kernels. As it was shown in \cite{Sc}, the natural framework for
studying hermitian kernels is given by Kre\u{\i}n spaces and for 
this reason we
briefly discuss the necessary terminology.

\subsection{Kre\u\i n Spaces}
An indefinite inner product space
$(\cH,[\cdot,\cdot]_\cH)$
is called {\it Kre\u{\i}n space}
provided that there exists 
a positive inner product $\langle\cdot,\cdot\rangle $
turning $(\cH, \langle\cdot,\cdot\rangle )$ into a Hilbert space
and such that $[\xi ,\eta ]=\langle J\xi ,\eta \rangle $, 
$\xi ,\eta \in \cH $,
for some symmetry $J$ ($J^*=J^{-1}=J$ with respect to the Hilbert space
structure) on $\cH$.
Such a symmetry $J$ is called a {\it fundamental symmetry} and we will
frequently indicate by a lower index the space on which it acts.
For two Kre\u{\i}n spaces $\cH$ and $\cK$ 
we denote by $\lhk$ the set 
of linear bounded operators from $\cH $ to $\cK$.
For $T\in \lhk$ we denote by $T^{\sharp }\in \lkh$ the adjoint 
of $T$ with respect to the indefinite inner product $[\cdot,\cdot]$.
The Hilbert space adjoint of $T$ with respect to the positive 
inner products $\langle\cdot,\cdot\rangle $ 
is denoted by $T^*$. It is important to note that, if $J_\cH$ and $J_\cK$ are
fundamentaly symmetries on $\cH$ and, respectively, $\cK$ then
\begin{equation*} T^\sharp =J_\cH T^* J_\cK.\end{equation*}

We say that $A\in \lh$ is a {\em selfadjoint 
operator} if $A^{\sharp }=A$. For example, in terms of fundamental symmetries,
this means $J_\cH A=A^* J_\cH$. Also, we say that the operator $U\in \lhk$ is 
{\it unitary} if $UU^{\sharp }=I_\cK$ and $U^{\sharp }U=I_\cH$, where 
$I_\cH$ denotes the identity operator on $\cH$. Equivalently, this means that
$U$ is boundedly invertible and $J_\cH U^{-1}=U^* J_\cK$. In terms of inner
products, this means that $U$ is isometric and surjective.

A special situation occurs for a unitary operator $U$ with domain and range the
same Kre\u\i n space $\cH$ if it commutes with some fundamental symmetry
$J_\cH$. Such a unitary operator is called {\em fundamentally reducible} and
it can be characterized in other different ways. For instance, $U$ is
fundamentally reducible if and only if it is power bounded.

Most of the difficulties in dealing with operators on 
Kre\u{\i}n spaces are caused by the 
lack of a well-behaved factorization 
theory. The concept of induced space
turned out to be quite useful in order to deal with this issue.
Thus, let $\cH$ be a Hilbert space and, for a selfadjoint operator $A$ 
in $\lh$,
we define a new inner product $[\cdot ,\cdot ]_A$
on $\cH$ by the formula
\begin{equation}\label{subA}
[\xi ,\eta ]_A=\langle A\xi ,\eta \rangle_\cH, \quad \xi ,\eta \in \cH .
\end{equation} 
A pair $(\cK, \Pi )$
consisting of a Kre\u{\i}n space $\cK$
and a bounded operator $\Pi \in \lhk$ 
is called a {\it Kre\u{\i}n space induced}
by $A$ provided that $\Pi $ has dense range
and the relation 
\begin{equation}\label{indus}
[\Pi \xi ,\Pi \eta ]_\cK=\langle \xi ,\eta \rangle_A 
\end{equation} 
holds for all $\xi ,\eta \in \cH$.
There are many known examples of induced spaces. 
A more delicate question is the {\em uniqueness} of the 
induced Kre\u\i
n spaces (see  \cite{CG1}).

\subsection{Hermitian Kernels}
We can use the concept of induced space in order
to describe the Kolmogorov decomposition of a 
hermitian kernel.
Let $X$ be an arbitrary set. From now on we assume $\cH$ 
is a Hilbert space with inner product denoted by 
$\langle\cdot,\cdot\rangle $. 
A {\it kernel } on $X$ is a mapping $K$ defined
on $X\times X$ with values in $\lh$. 
The adjoint $K^*$ of $K$ is defined by the formula
$K^*(x,y)=K(y,x)^*$.
The kernel $K$ is called {\it hermitian on} 
$X$ if
$K^*=K$. 

Let ${\cF}_0(X,\cH)$ denote the vector space of all functions
on $X$ with values in $\cH$ which vanish except on a 
finite number of points. 
We associate to $K$ an inner product on 
${\cF}_0(X,\cH)$
by the formula:
\begin{equation}\label{fega}
[f,g]_K=\sum_{x,y\in X}
\langle K(x,y)f(y),g(x)\rangle,\quad f,g\in{\cF}_0(X,\cH). 
\end{equation}
We say that the hermitian kernel $L:X\times X\rightarrow \lh$
is {\it positive definite} if the inner product
$[\cdot ,\cdot ]_L$ associated to $L$ 
by the formula \eqref{fega} is positive. 
One can introduce a natural partial order on 
the set of hermitian kernels on $X$ with values in $\lh$
as follows:
if $A$, $B$ are hermitian kernels, then $A\leq B$ means
$[f,f]_A\leq [f,f]_B$ for all $f\in {\cF}_0(X,\cH)$.

A {\it Kolmogorov decomposition}
of the hermitian kernel $K$ is a pair $(V;\cK)$ with the following properties:
\begin{itemize}
\item[KD1] $\cK$ is a Kre\u{\i}n space with fundamental symmetry $J$;
\item[KD2] $V=\{V(x)\}_{x\in X}\subset\lhk $ such that
$K(x,y)=V(x)^{*}JV(y)$ for all $x,y\in\cK$;
\item[KD3] $\{V(x){\cH}\mid x\in X\}$ is total in $\cK$.
\end{itemize}

The next result, obtained in \cite{CG1}, settles the question
concerning the existence of a Kolmogorov decomposition 
for a given hermitian kernel.

\begin{theorem}\label{kolmo}
Let $K:X\times X\rightarrow \lh$
be a hermitian kernel. The following assertions are equivalent:

{\rm (1)} There exists 
a positive definite kernel $L:X\times X\rightarrow \lh$
such that $-L\leq K\leq L$.

{\rm (2)} $K$ has a Kolmogorov decomposition.
\end{theorem}

The condition in assertion {\rm (1)} of the 
previous result appeared earlier in the work of 
L.~Schwartz \cite{Sc} concerning the structure of hermitian
kernels. It is easy to see that {\rm (1)}
is also equivalent to the representation of $K$ as a difference of two 
positive definite kernels. Thus, Theorem~\ref{kolmo} says that the class of
hermitian kernels admitting Kolmogorov decompostions is the same with the
class of hermitian kernels in the linear span of the cone of positive definite
kernels. 

It is convenient for our purpose to review 
a construction of Kolmogorov decompositions.
We assume that there exists a positive definite kernel 
$L:X\times X\rightarrow \lh$ such that $-L\leq K\leq L$.
Let ${\cH}_L$ be the Hilbert space obtained by the completion of 
the quotient space
${\cF}_0(X,\cH)/{\cN}_L$ with respect to $[\cdot ,\cdot ]_L$, 
where ${\cN}_L=\{f\in {\cF}_0(X,\cH)\mid [f,f]_L=0\}$
is the isotropic subspace of the inner product space
$({\cF}_0(X,\cH),[\cdot ,\cdot ]_L)$.
Since $(1)$ in Theorem~\ref{kolmo} is equivalent to
\begin{equation*}\label{schwarz}
|[f,g]_K|\leq [f,f]_L^{1/2}\, [g,g]_L^{1/2}
\end{equation*}
for all $f,g\in{\cF}_0(X,\cH)$ (see Proposition~38, \cite{Sc}), 
it follows that ${\cN}_L$ is a subset of the 
isotropic subspace ${\cN}_K$ of the inner product space
$({\cF}_0(X,{\cH}),[\cdot ,\cdot ]_K)$.
Therefore, $[\cdot ,\cdot ]_K$ uniquely induces an inner product on 
${\cH}_L$, still denoted by $[\cdot ,\cdot ]_K$, such that
\eqref{schwarz} holds for $f,g\in{\cH}_L$.
By the Riesz representation theorem we obtain a 
selfadjoint contractive operator $A_L\in{\cL}({\cH}_L)$,
referred to as the {\it Gram operator} of $K$ with respect to $L$, 
such that 
\begin{equation*}\label{gram}
[f,g]_K=[A_Lf,g]_L,\quad f,g\in{\cF}_0(X,\cH).
\end{equation*}
Let $(\cK ,\Pi )$ be a
Kre\u{\i}n space induced by $A_L$.
For $\xi \in \cH$ and $x\in X$, 
we define the element $\xi _x=\delta_x \xi\in {\cF}_0(X,\cH)$ (here $\delta_x$
is the Kronecker function delta), that is,
\begin{equation}\label{xi}\xi_x(y)=\left\{\begin{array}{ll}
\xi, & y=x; \\
0, & y\ne x.
\end{array}\right.\end{equation}
Then we define
\begin{equation*}\label{v}
V(x)\xi =\Pi [\xi _x],\quad x\in X,\ \xi\in\cH,
\end{equation*}
where $[\xi _x]=\xi_x+\cN_L$ denotes the class of $\xi _x$ in ${\cH}_L$ and 
it can be verified that $(V;\cK)$ is a Kolmogorov decomposition of 
the kernel $K$.


We finally review 
the uniqueness property of the Kolmogorov decomposition.
Two Kolmogorov decompositions $(V_1,\cK _1)$ and
$(V_2,\cK _2)$ of the same hermitian
kernel $K$ are {\it unitarily
equivalent} if there exists a unitary operator 
$\Phi\in\cL(\cK _1,\cK _2)$
such that for all $x\in X$ we have $V_2(x)=\Phi V_1(x)$.
The following result was obtained in \cite{CG1}. We denote by $\rho(T)$ the
resolvent set of the operator $T$.

\begin{theorem}\label{unic}
Let $K$ be a hermitian kernel which has Kolmogorov
decompositions. The following assertions are equivalent:\smallskip

{\rm (1)} All Kolmogorov decompositions of $K$ are unitarily equivalent.

{\rm (2)} For each positive definite
kernel $L$ such that $-L\leq K\leq L$, there exists $\epsilon>0$ such that
either $(0,\epsilon)\subset\rho(A_L)$ or $(-\epsilon,0)\subset\rho(A_L)$,
where $A_L$ is the Gram operator of $K$ with respect to $L$.
\end{theorem}

\subsection{Invariant Hermitian Kernels}
We now review some results on the Kolmogorov decomposition of
hermitian kernels with additional symmetries.
Let $\phi$ be an action of a unital semigroup $S$ on $X$.
Assume that the $\lh$-valued hermitian kernel $K$ has a Kolmogorov
decomposition
$(V,\cK)$. The action $\phi $ is linearized by the following mapping:
for any $a\in S$, $x\in X$ and $\xi \in \cH$,
\begin{equation}\label{rep}
U(a)V(x)\xi =V(\phi (a,x))\xi.
\end{equation}
We notice that for $a,b\in S$, $x\in X$ and $\xi\in\cH$ we have
\begin{equation*}\begin{split}
U(a)U(b)V(x)\xi &=U(a)V(\phi (b,x))\xi =V(\phi (a,\phi (b,x)))\xi  \\
 &=V(\phi (ab,x))\xi =U(ab)V(x)\xi.
\end{split}\end{equation*}
Therefore, the family $\{U(a)\}_{a\in S}$ is a semigroup of linear 
operators with a common dense domain $\bigvee _{x\in X}V(x)\cH$
(throughout this paper $\bigvee $ denotes the linear
space generated by some set, without taking any closure).
If $K$ is a positive definite kernel then the previous construction
is well-known (see, for instance, \cite{PS}). The remaining question,
especially in case $K$ is not positive definite,
is: what additional conditions on the kernel $K$ should be imposed
in order to ensure the boundedness of the operators $U(a)$, $a\in S$? 
We gave a possible answer in \cite{CG2}, by considering an additional symmetry
of the kernel.

Consider the set $B=\{\xi _x\mid \xi \in {\cH}, x\in X\}$ 
which is a vector space basis of ${\cF}_0(X,\cH)$.
Define for $a\in S$,
\begin{equation}\label{actpsi}
\psi _a(\xi _x)=\xi _{\phi (a,x)}
\end{equation}
and this mapping can be extended by linearity to a linear
mapping, also denoted by $\psi _a$, from ${\cF}_0(X,\cH)$
into ${\cF}_0(X,\cH)$.
We say that a positive definite kernel $L$ is 
$\phi $-{\it bounded} provided that for all $a\in S$, 
$\psi _a$ is bounded with respect to the seminorm
$[\cdot ,\cdot ]_L^{1/2}$ induced by $L$ on ${\cF}_0(X,\cH)$.
We denote by 
${\cB}_{\phi }^+(X,\cH)$ the set of positive definite $\phi$-bounded
kernels on $X$ with values in $\lh$.

From now on we assume that $S$ is a unital semigroup 
with involution, that is, there exists a 
mapping $\fI:S\rightarrow S$ such that $\fI ^2=\,$the 
identity on $S$, and $\fI(ab)=\fI(b)\fI(a)$ for 
all $a,b\in S$. 
The following result was obtained in \cite{CG2}.
\begin{theorem}\label{baza}
Let $\phi $ be an action of the unital semigroup 
$S$ with involution $\fI$ on the set $X$ and let $K$ be an $\lh$-valued 
hermitian kernel on $X$ with the property that 
\begin{equation}\label{simetria}
K(x,\phi (a,y))=K(\phi (\fI (a),x),y)
\end{equation}
for all $x,y\in X$ and $a\in S$.
The
following assertions are equivalent:
\smallskip

{\rm (1)} There exists $L\in {\cB}_{\phi }^+(X,\cH)$ such that $-L\leq K\leq
L$.\smallskip 

{\rm (2)} $K$ has a Kolmogorov decomposition $(V;\cK)$
with the property that there exists a 
representation $U$ of $S$ on $\cK$ 
such
that 
\begin{equation}\label{rel}
V(\phi (a,x))=U(a)V(x)
\end{equation}
for all $x\in X$, $a\in S$.
In addition, $U(\fI (a))=
U(a)^{\sharp }$ for all $a\in S$. \smallskip 

{\rm (3)} $K=K_1-K_2$ for two positive definite 
kernels such that $K_1+K_2\in {\cB}^+_\phi(X,\cH)$.\smallskip 

{\rm (4)} $K=K_+-K_-$ for two disjoint positive definite 
kernels such that $K_++K_-\in {\cB}^+_\phi(X,\cH)$.
\end{theorem}

\subsection{Reproducing kernel spaces}
We now describe another construction, closely related to the 
Kolmogorov decomposition of a hermitian kernel.
Let $K$ be a hermitian kernel satisfying \eqref{simetria}
with a Kolmogorov decomposition $(V;\cK)$ 
as in Theorem~\ref{baza}. Define
$${\cR}=
\{g_f:X\rightarrow {\cH}
\mid g_f(x)=V^{\sharp }(x)f, \, f\in {\cK}\}.$$
Then $\cR$ is a vector subspace
of ${\cF}(X, \cH)$, the class of functions defined 
on $X$ with values in $\cH$. We define a map 
$\Phi :{\cK}\rightarrow \cR$ by 
$$\Phi f=g_f,\quad f\in\cK.$$
This map is linear and bijective, so that we can define
on $\cR$ the inner product
$$[g_f,g_h]_{\cR}=[f,h]_{{\cK}},\quad f,h\in\cK.$$
One checks that $\cR$ is a Kre\u{\i}n space with respect
to this inner product. Also, $\Phi $ is a bounded
 operator between the Kre\u{\i}n spaces 
$\cK$ and $\cR$, since it is closed and everywhere defined on 
$\cK$, hence it is unitary. Moreover, $\cR$
is the closure of the linear space generated by the 
functions $g_{V(y)\xi }$, $y\in X$ and $\xi \in \cH$.
These functions are related to the kernel $K$
as follows:
$$g_{V(y)\xi }(x)=V^{\sharp }(x)V(y)\xi =K(x,y)\xi,\quad x,y\in X,\, 
\xi\in\cH.$$
We will write $g_{y,\xi }$ instead of 
$g_{V(y)\xi }$, and since $g_{y,\xi }(x)=K(x,y)\xi $, 
these functions can be defined without using $V$. 
Therefore, the space $\cR$ has the following
{\it reproducing property}:
\begin{equation}\label{repr}
[g_f(x),\xi ]_{\cH}=[g_f,g_{x,\xi }]_{\cR},\quad x\in X,\ f\in\cK,\ \xi\in\cH.
\end{equation}
We also note that property \eqref{simetria}
of the kernel $K$ 
is reflected into a certain symmetry of the elements of $\cR $.
Thus, we define an operator $\bar U(a)\in {\cL}(\cR)$ by
$$\bar U(a)=\Phi U(a)\Phi ^{\sharp },\quad a\in X,$$
where $U$ is the projective representation of $S$ given
by Theorem~\ref{baza}.
We have 
$$\bar U(a)g_f=\Phi U(a)\Phi ^{\sharp }g_f=\Phi U(a)f=g_{U(a)f},\quad a\in X,\
f\in\cK.$$ 
On the other hand, for any $a,x\in X$ and $f\in\cK$,
\begin{equation*}\begin{split}
g_f(\phi (\fI (a),x))&=V^{\sharp }(\phi (\fI (a),x))f=
V(x)^{\sharp }U(\fI (a))^{\sharp }f \\
&=g_{U(a)^{\sharp }f}(x)
\end{split}\end{equation*}
and we deduce that the elements of $\cR$ satisfy the relation
$$
(\bar U(a)g_f)(x)=\sigma (\fI (a),a)^{-1}
g_f(\phi (\fI (a),x)).
$$
Based on this relation we obtain the following result.
\begin{theorem}\label{repro}
Let $\phi $ be an action of the unital semigroup 
$S$ with involution $\fI$ on the set $X$ and let $K$ be an $\lh$-valued 
hermitian kernel on $X$ with the property that 
\begin{equation}\label{sim}
K(x,\phi (a,y))=K(\phi (\fI (a),x),y)
\end{equation}
for all $x,y\in X$ and $a\in S$.
The
following assertions are equivalent:
\smallskip

{\rm (1)} There exists $L\in {\cB}_{\phi }^+(X,\cH)$ such that $-L\leq K\leq
L$.\smallskip

{\rm (2)} $K$ has a Kolmogorov decomposition $(V;\cK)$
with the property that there exists a 
representation $U$ of $S$ on $\cK$ 
such
that $V(\phi (a,x))=U(a)V(x)$
for all $x\in X$, $a\in S$.\smallskip

{\rm (3)} There exists a Kre\u{\i}n space $\cR$ such that

\quad $(a)\quad {\cR}\subset {\cF}(X,{\cH}).$

\quad $(b)$\quad 
The set $\{g_{x,\xi}\mid x\in X, \xi \in \cH\}$ is total in 
${\cR}$.

\quad $(c)\quad [f(x),\xi ]_{\cH}=[f,g_{x,\xi }]_{\cR}$ for all
$f\in {\cR}, \xi \in {\cH}, x\in X.$

\quad $(d)$\quad There exists a representation
$\bar U$ of $S$ on $\cR$ such that 
$$(\bar U(a)f)(x)=
f(\phi (\fI (a),x)) $$
for all $a\in S$, $x\in X$ and $f\in \cR$.  
\end{theorem}
\begin{proof} The implication $(2)\Rightarrow (3)$ was already 
proved above. In order to prove $(3)\Rightarrow (2)$
we define the linear mapping $\bar V(x)$ from $\cH$ 
into $\cR$ by the formula:
\begin{equation*}\bar V(x)\xi =g_{x,\xi },\quad x\in X,\ 
\xi\in\cH.\end{equation*}
The property $(c)$ shows that $\bar V(x)$ is a closed operator and by the 
closed graph theorem we deduce that 
$\bar V(x)\in {\cL}(\cH, \cR)$.
From $(b)$ and $(c)$ we deduce that $(\bar V;{\cR})$ 
is a Kolmogorov decomposition of $K$. Finally, 
\begin{equation*}\begin{split}
(\bar V(\phi (a,x))\xi )(y)&=g_{\phi (a,x),\xi }(y)=
K(y,\phi (a,x))\xi  \\
&=K(\phi (\fI (a),y),x)\xi  \\
& =g_{x,\xi }(\phi (\fI (a),y))\\
&= (\bar U(a)\bar V(x)\xi )(y). 
\end{split}\end{equation*}
This completes the proof. 
\end{proof}

\section{Hankel type kernels}
In this section we interpret the invariance property 
\eqref{simetria} as a Hankel condition.
To see this, let $S=\NN $ be
the additive semigroup of natural numbers (including $0$) 
and the action $\phi $
is given by right translation. If $K$ satisfies
\eqref{simetria}, then
$K(n,p+m)=K(p+n,m)$ for $m,n,p\in \NN$ and $K$ is a 
so-called {\em Hankel kernel}. We can extend this example to a noncommutative
setting as follows.
Let $S=\FF _N^+$ be the unital free semigroup 
on $N$ generators $g_1,\ldots ,g_N$ with lexicograhpic
order $\prec $. The empty word is the identity element
and the length of the word $\sigma $ is 
denoted by $|\sigma |$. The length of the empty word is $0$.
There is a natural involution on
$\FF _N^+$
given by $\fI(g_1\ldots g_k)=g_k\ldots g_1$ as well as a natural
action of $\FF _N^+$ on itself by juxtaposition, $\phi (\sigma,\tau )=
\sigma \tau $, $\sigma ,\tau \in \FF _N^+$.
The condition  
\eqref{simetria} means in this case that 
\begin{equation}\label{hankel}
K(\sigma ,\beta \tau )=K(\fI(\beta )\sigma ,\tau)
\end{equation}
for $\beta ,\sigma ,\tau \in \FF _N^+$.
It was noticed in \cite{Co} that kernels as
above appear in connection with orthogonal 
polynomials in $N$ indeterminates satisfyting the 
relations $Y^*_k=Y_k$, $k=1,\ldots ,N$.

Let $\cP_N^0$ be the algebra of polynomials
in $N$ non-commuting
indeterminates $Y_1$,$\ldots $,$Y_N$ 
with complex coefficients. For any 
$\sigma=g_{j_1}g_{j_2}\cdots g_{j_l}\in \FF_N^+$, 
where $j_p\in\{1,2,\ldots,N\}$ for all $p=1,\ldots,l$, $l=|\sigma|$, we denote
$Y_\sigma=Y_{g_{j_1}} Y_{g_{j_2}} \cdots Y_{g_{j_l}}$. With this notation, 
each element $P\in \cP_N^0$ can be uniquely written 
as
\begin{equation}\label{pesi}
P=\sum _{\sigma \in \FF _{N}^+}c_{\sigma }Y_{\sigma },
\end{equation} where $(c_{\sigma })_{\sigma\in\FF_N^+}\subset \CC$ 
has finite support.

An involution $*$ 
on $\cP _N^0$ can be introduced as follows: $Y_k^*=Y_{k}$, $k=1,\ldots ,N$;
on monomials, 
$(Y_{\sigma })^*=Y_{\fI(\sigma )}$; and, in general,
 if $P$ has the representation as in \eqref{pesi}
then 
\begin{equation*}P^*=\sum _{\sigma \in \FF _{N}^+}\overline{c}_{\sigma }
Y_{\sigma }^*.\end{equation*}
Thus, $\cP_N^0$ is a unital, associative, $*$-algebra over $\CC $.
A linear functional $Z$ on $\cP_N^0$ is called {\em hermitian}
if $Z(P^*)=\overline{Z(P)}$ for $P\in \cP_N^0$.

A convenient subclass of hermitian functionals, 
called {\em GNS functionals}, is given by those functionals admitting
GNS data. A triplet $(\pi ,\cK, \Omega )$ is called
a {\em GNS data} associated to $Z$ 
if $\pi $ is a hermitian closable representation
of $\cP_N^0$ on a Kre\u{\i}n space $\cK$
and $\Omega \in \cD (\pi )$, the domain of $\pi $, 
such that $Z(P)=[\pi (P)\Omega ,\Omega ]_{\cK }$
for $P\in \cP_N^0$, and $\bigvee _{P\in \cP_N^0}\pi (P)\Omega 
=\cD (\pi )$ (see \cite{AGW}, \cite{Ho}). 
The numbers
$s_{\sigma }=Z(Y_{\sigma })$, $\sigma \in \FF _{N}^+$,
are called the {\em moments} of $Z$.

Conversely, to any family of complex numbers
 $\Sigma=(s_\sigma)_{\sigma\in\FF_N^+}$,
we can associate the kernel 
\begin{equation}\label{kaces}K_\Sigma(\sigma ,\tau )=
s_{\fI(\sigma )\tau },\quad\sigma ,\tau \in \FF _{N}^+,\end{equation}
and it is easy to see that this kernel  
satisfies \eqref{hankel}.

The following is a Hamburger type description of moments.

\begin{theorem}\label{hamburger}
The complex numbers $s_{\sigma }$, $\sigma \in \FF _{N}^+$,
are the moments of a GNS functional on $\cP_N^0$
if and only if there exists a positive definite kernel
$L$ on $\FF _{N}^+$
such that $-L\leq K\leq L$, where $K(\sigma ,\tau )=
s_{\fI(\sigma )\tau }$, $\sigma ,\tau \in \FF _{N}^+$.
\end{theorem}

\begin{proof}
This result is just another facet of Theorem~\ref{kolmo}.
Assume first that the numbers
$s_{\sigma }$, $\sigma \in \FF _{N}^+$,
are the moments of a GNS functional on $\cP_N^0$.
Let $(\pi ,\cK, \Omega )$ be a GNS data associated to $Z$.
Define $V:\FF _{N}^+\rightarrow \cL (\CC ,\cK )$ by the formula:
$$V(\sigma )\lambda =\pi (Y_{\sigma })(\lambda \Omega ),\quad 
\sigma \in \FF _{N}^+, \lambda \in \CC.$$
We deduce that for  $\sigma ,\tau \in \FF _{N}^+$
and  $\lambda ,\mu \in \CC,$
\begin{equation*}\begin{split}
V(\sigma )^{\sharp }V(\tau )\lambda \overline{\nu }&=
[V(\tau )\lambda ,V(\sigma )\nu ]_{\cK } 
 =[\pi (Y_{\tau })(\lambda \Omega ),
\pi (Y_{\sigma })(\nu \Omega )]_{\cK } \\
 &=\lambda \overline{\nu }
[\pi (Y_{\sigma }^*Y_{\tau })\Omega ,\Omega ]_{\cK }
 =\lambda \overline{\nu }Z(Y_{I(\sigma )\tau })=
\lambda \overline{\nu }K(\sigma ,\tau ).
\end{split}
\end{equation*}
Also, the set $\{V(\sigma )\lambda \mid 
\sigma \in \FF _{N}^+, \lambda \in \CC\}$ 
is total in $\cK$, so that $(V,\cK)$ is a Kolmogorov
decomposition of $K$. By Theorem~\ref{kolmo}, 
there exists a positive definite kernel
$L$ on $\FF _{N}^+$
such that $-L\leq K\leq L$.

Conversely, let $(V,\cK)$ be a Kolmogorov
decomposition of $K$. Define $\Omega =V(\emptyset )$ and 
$$\pi (Y_{\sigma })\Omega =V(\sigma ), \quad \sigma \in \FF _{N}^+.$$
We notice that $\bigvee _{P\in \cP_N^0}\pi (P)\Omega
=\bigvee _{\sigma \in \FF _{N}^+}V(\sigma )\CC $, we define
$\cD (\pi )=\bigvee _{\sigma \in \FF _{N}^+}V(\sigma )\CC $,
and we can extend $\pi $ to $\cP_N^0$
by linearity. Clearly, $\cD (\pi )$
is invariant under $\pi (P)$, $P\in \cP_N^0$, 
and $\pi (P)\pi (Q)=\pi (PQ)$. Also, for $k,k'\in \cD (\pi )$,
\begin{equation*}\begin{split}
[\pi (Y_{\sigma })k,k']_{\cK }&=
[\pi (Y_{\sigma })\sum _{k=1}^nc_k\pi (Y_{\tau _k})\Omega ,
\sum _{j=1}^md_j\pi (Y_{\tau '_j})\Omega ]_{\cK }\\
&=\sum _{k,j=1}^{n,m}c_k\overline{d}_j
[\pi (Y_{\sigma })\pi (Y_{\tau _k})\Omega ,\pi (Y_{\tau '_j})\Omega ]_{\cK }\\
&=\sum _{k,j=1}^{n,m}c_k\overline{d}_j[V(\sigma \tau _k)1,
V(\tau '_j)1]_{\cK }
=\sum _{k,j=1}^{n,m}c_k\overline{d}_jV(\tau '_j)^{\sharp }V(\sigma \tau _ k)\\
&=\sum _{k,j=1}^{n,m}c_k\overline{d}_jK(\tau '_j,\sigma \tau _k)= 
\sum _{k,j=1}^{n,m}c_k\overline{d}_jK(\fI(\sigma )\tau '_j,\tau _k) \\
 &=[k,\pi (Y_{\fI(\sigma )})k']_{\cK },
\end{split}\end{equation*}
which shows that the domain of $\pi (Y_{\sigma })^{\sharp }$
contains $\cD (\pi )$ and 
$$\pi (Y_{\sigma })^{\sharp }|\cD (\pi )=\pi (Y_{\fI(\sigma )})
=\pi (Y_{\sigma }^*).$$
We can extend this argument and show that the same is true for any 
$P\in \cP_N^0$, so that $(\pi ,\cK ,\Omega )$ is a GNS data
for $Z(P)=[\pi (P)\Omega ,\Omega ]_{\cK }$, $P\in \cP_N^0$.
The moments of $Z$ are
\begin{equation*}\begin{split}
Z(Y_{\sigma })&=[\pi (Y_{\sigma })\Omega ,\Omega ]_{\cK }
 =[V(\sigma )1,V(\emptyset )1]_{\cK } \\
 &= V(\emptyset )^{\sharp }V(\sigma )=K(\emptyset ,\sigma )
=s_{I(\emptyset )\sigma }=s_{\sigma }.
\end{split}\end{equation*}
\end{proof}

As a consequence of the previous result and Theorem~\ref{unic}, 
we deduce a uniqueness condition for the solvability
of the Hamburger moment problem for GNS functionals.

\begin{theorem}\label{ham}
Let $s_{\sigma }$, $\sigma \in \FF _{N}^+$,
be the moments of some GNS
functional on  $\cP_N^0$, and consider the kernel $K(\sigma ,\tau )=
s_{\fI(\sigma )\tau }$, $\sigma ,\tau \in \FF _{N}^+$.
Then, there exists a unique GNS functional on $\cP_N^0$
with moments $s_{\sigma }$, if and only if for each positive
definite kernel $L$ on $\FF_{N}^+$
such that $-L\leq K\leq L$, 
there exists $\epsilon >0$ such that either
$(0,\epsilon )\subset \rho (A_L)$ or
$(-\epsilon ,0)\subset \rho (A_L)$, where $A_L$ is the Gram operator
of $K$ with respect to $L$.
\end{theorem}

\section{Dilations and decomposition of kernels}
In this section we show that Theorem~\ref{baza} provides a general framework
for a version of the Stinespring theorem and for decompositions 
of hermitian linear maps.
Let $\cA$ be a unital $*$-algebra, $\cH$ a Hilbert space, 
and let $T:\cA\rightarrow \lh$ be a linear
hermitian map. A {\it Stinespring
dilation} of $T$ is, by definition, a triplet $(\pi,\cK,B)$ such that:

\begin{itemize}
\item[SD1] $\cK$ is a Kre\u\i n space with a fundamental symmetry $J$, and
$B\in\cL(\cH,\cK)$;
\item[SD2] $\pi\colon\cA\rightarrow \cL(\cK)$ is a
selfadjoint (that is, $\pi(a^*)=\pi(a)^\sharp=J\pi(a)^*J$ 
for all $a\in\cA$) representation,  such that $
T(a)=B^{*}J\pi (a)B$, for all $a\in\cA$.
\end{itemize}
If, in addition, 
\begin{itemize}
\item[SD3] $\{\pi(a)B\cH \mid a\in \cA\}$ is total in $\cK$,
\end{itemize} 
then the Stinespring dilation is called {\it minimal}.

We consider the action $\phi$ of $\cA$ on itself defined by 
\begin{equation}\label{fi}
\phi (a,x)=xa^*,\quad x,a\in \cA,
\end{equation} 
and a hermitian kernel is associated to $T$ by the formula
\begin{equation}\label{kate}K_T(x,y)=T(xy^*),\quad x,y\in\cA.\end{equation}
It readily follows that $K_T$ is $\phi$-invariant, that is,
\begin{equation}\label{siminica}
K_T(x,\phi (a,y))=T(xay^*)=K_T(\phi (a^*,x),y),\quad a,x,y\in\cA.
\end{equation}

\begin{proposition}\label{stinekol} Given a minimal Stinespring dilation 
$(\pi,\cK,B)$ of the hermitian linear map $T\colon\cA\rightarrow\cL(\cH)$,
let
\begin{equation}\label{vepi}V(x)=\pi(x^*)JB,\quad x\in\cA,\end{equation}
where $J$ is a fundamental decomposition of $\cK$.
Then $(V,\cK)$ is an invariant Kolmogorov decomposition of the
hermitian kernel $K_T$. 

In addition, \eqref{vepi} establishes a
bijective correspondence between the set of minimal Stinespring
dilations of $T$ and the set of invariant Kolmogorov decompositions of
$K_T$.\end{proposition}
 
\begin{proof} Let $(\pi,\cK,B)$ be a minimal Stinespring dilation of
$T$ and define $(V,\cK)$ as in \eqref{vepi}. Then
$$V(x)^*J V(y)=B^*J \pi(x^*)^*J \pi(y^*)JB=B^*J
\pi(xy^*)JB=T(xy^*)=K_T(x,y),\quad x,y\in\cA.$$
Since $\bigvee_{a\in\cA}\pi(a)B\cH=\bigvee_{x\in\cA}V(x)\cH$ it
follows that $(V,\cK)$ is a Kolmogorov decomposition of $T$.
Let us note that, by the definition of $V$,
$$\pi(a)V(x)=\pi(a)\pi(x^*)JB=\pi(ax^*)JB=V(xa^*)\,\quad a,x\in\cA,$$
and hence, letting $U=\pi$, it follows that the Kolmogorov
decomposition $(V,\cK)$ is $\phi$-invariant.

Conversely, let $(V,\cK)$ be an invariant Kolmogorov decomposition of the
hermitian kernel $K_T$, that is, there exists a hermitian
representation $U\colon\cA\rightarrow\cL(\cK)$ 
of the multiplicative semigroup $\cA$ with involution $*$, such that
$$U(a)V(x)=V(xa^*),\quad a,x\in\cA.$$
Define $\pi=U$ and $B=JV(1)$. Since $T$ is linear it follows easily
that $\pi$ is also linear, hence a selfadjoint representation of the
$*$-algebra $\cA$ on the Kre\u\i n space $\cK$. Then, taking into account
that $V(a)=U(a^*)B$ for all $a\in\cA$, it follows
$$T(a)=V(a)^*J V(1)=B^* U(a^*)^* JB=B^* JU(a)B,\quad
a\in\cA,$$ and since
$\bigvee_{a\in\cA}\pi(a)B\cH=\bigvee_{x\in\cA}V(x)\cH$ we thus proved
that $(\pi,\cK,B)$ is a minimal Stinespring dilation of $T$.
One easily check that the mapping defined in
\eqref{vepi} is the inverse of the mapping associating to each
invariant Kolmogorov decomposition $(V,\cK)$ the minimal Stinespring
dilation $(\pi,\cK,B)$ as above.\end{proof}

Proposition~\ref{stinekol} reduces the existence of Stinespring dilations of
hermitian maps $T$ to the existence of invariant Kolmogorov decompositions for
the hermitian kernel $K_T$ defined as in \eqref{kate}. Now the following 
result is just an application of 
Theorem~\ref{baza}.

\begin{theorem}\label{stine}
Let $\cA$ be a unital $*$-algebra  
and let $T:\cA\rightarrow \lh$ 
be a linear hermitian map. The following assertions are equivalent:
\smallskip

{\rm (1)} There exists a positive definite kernel 
$L\in \cB _{\phi }^+(\cA,\cH)$, $\phi $ given by \eqref{fi}, 
such that $-L\leq K_T\leq L$.

{\rm (2)} $T$ has a minimal Stinespring dilation. 
\end{theorem}

We show now that Wittstock's Decomposition
Theorem \cite{Wi} and Paulsen's Dilation Theorem \cite{Pa} fit into the
framework of invariant Kolmogorov decompositions, more precisely, 
the following result shows that in case $\cA$ is a $C^*$-algebra,
Schwartz's boundedness condition for hermitian kernels represents
an extension of the concept of completely bounded map.
We use standard terminology from the theory of operator spaces, 
e.g.\ see \cite{Pa} and \cite{ER}.

\begin{theorem}\label{witt}
Let $\cA$ be a unital $C^*$-algebra and let $T:\cA \rightarrow \lh$ 
be a linear hermitian map. The following assertions are equivalent:

{\rm (1)} $T$ is completely bounded.\smallskip

{\rm (2)} There exists a completely positive map $S:\cA \rightarrow \lh$
such that $-S\leq T\leq S$.\smallskip

{\rm (3)} There exists a Hilbert space $\cK$ with a symmetry $J$, 
a $*$-representation 
$\pi :\cA\rightarrow \lk $ commuting with $J$, 
and a bounded operator $B\in \lhk $
such that 
\begin{equation*}T(a)=B^*J\pi (a)B,\quad a\in \cA,\end{equation*} and 
$\bigvee _{a\in \cA}\pi (a)B\cH$ is dense in $\cK$.\smallskip
 
{\rm (4)} $T=T_+-T_-$ for two completely positive maps $T_+$ and 
$T_-$.
\end{theorem}
\begin{proof}
In the following we let $\cU(\cA)$ be the unitary group of $\cA$.
Then $\cU(\cA)$ has the involution $\fI (a)=a^{-1}=a^*$
and acts on $\cA$ by $\phi (a,x)=xa^*=xa^{-1}$.

$(1)\Rightarrow (2)$ 
 We use 
Paulsen's off-diagonal technique. Briefly, 
assume that $T$ is completely bounded. By Theorem~7.3 in 
\cite{Pa}, there exist completely positive maps
$\phi _1$ and $\phi _2$ such that the map
$$F(\left[\begin{array}{cc} a & b\\ c &
d
\end{array}\right])=
\left[\begin{array}{cc}
 \phi _1(a) & T(b) \\
 T(c^*)^* & \phi _2(d) 
\end{array}\right]$$
is completely positive. Define $S(a)=\frac{1}{2}(\phi_1(a)+\phi_2(a))$, 
which is a completely positive map. We can check that 
$-S\leq T\leq S$. First, let $a\geq 0$, $a\in \cA$.
Then 
$\left[\begin{array}{cc}
a & \pm a \\
 \pm a & a 
\end{array}\right]\geq 0$,
so that
$\left[\begin{array}{cc}
\phi _1(a) & \pm T(a)\\
\pm T(a) & \phi _2(a)
\end{array}\right]\geq 0$. In particular, for $\xi \in \cH$,
$$\langle \left[\begin{array}{cc}
\phi _1(a) & \pm T(a)\\
\pm T(a) & \phi _2(a)
\end{array}\right]
\left[\begin{array}{cc} 
\xi \\ \xi \end{array}\right],
\left[\begin{array}{cc} 
\xi \\ \xi \end{array}\right]
\rangle \geq 0,$$
equivalently, $\langle (\phi _1(a)\pm 2T(a)+\phi _2(a))\xi ,\xi\rangle \geq 0$.
Therefore, $S\pm T$ are positive maps. The argument can
be extended in a straightforward manner (using the so-called
canonical shuffle as in \cite{Pa}) to show that 
$S\pm T$ are completely positive maps.

$(2)\Rightarrow (3)$ Since $S$ is completely positive, the 
kernel $K_S$ is positive definite and satisfies 
$-K_S\leq K_T\leq K_S$. Also, 
$$K_S(x,\phi (a,y))=K_S(\phi (a^{-1},x),y),\quad a\in\cU(\cA),\ x,y\in\cA.$$ 
By Theorem~4.3 in \cite{CG2}, there exists a 
Kolmogorov decomposition $(V,\cK )$
of $K_T$ and a fundamentally reducible representation
$U$ of $\cU(\cA )$ on $\cK $ such that
$$U(a)V(x)=V(xa^{-1}),\quad a\in\cU(\cA),\ x\in\cA.$$
 Let $J$ be a fundamental 
symmetry on $\cK$ such that $U(a)J=JU(a)$ for all $a\in \cU(\cA)$.
Then $U$ is also a representation of $\cU(\cA )$ on the Hilbert
space $(\cK ,\langle \cdot ,\cdot \rangle _J)$. Also, for $a\in \cU(\cA)$,
$$T(a)=K_T(a,1)=V(a)^{*}JV(1)=V(1)^*JU(a)V(1).$$
Since $\cA$ is the linear span of $\cU(\cA )$ and $T$ is linear, 
$U$ can be extended by linearity to a representation $\pi $ 
of $\cA$ on $\cK$ commuting with $J$ and such that  
$$T(a)=V(1)^*J\pi (a)V(1)$$
holds for all $a\in \cA$. Also, $\bigvee _{a\in \cU(\cA )}U(a)V(1)\cH
=\bigvee _{a\in \cU(\cA )}V(a^{-1})\cH =\bigvee _{a\in \cU(\cA )}V(a)\cH $
and using once again the fact that $\cA$ is the linear span of $\cU(\cA)$, 
we deduce that $\bigvee_{a\in \cA }U(a)V(1)\cH $ is dense in $\cK$. 
 
(3)$\Rightarrow$(4)
We define for $a\in \cA $,
$$V(a)=\pi (a^*)B.$$
Then $V(a)$ is in $\lhk $ and 
$\bigvee _{a\in \cA}V(a)\cH =\bigvee _{a\in \cA}\pi (a)B\cH $.
$\cK $ becomes a Kre\u{\i}n space by setting $[x,y]_{\cK }=
\langle Jx,y\rangle $, $x,y\in \cK $.
Also, for $\xi ,\eta \in \cH $,
\begin{equation*}\begin{split}
\langle V(x)^{*}JV(y)\xi ,\eta \rangle =&
\langle JV(y)\xi ,V(x)\eta \rangle 
= \langle J\pi (y^*)B\xi ,\pi (x^*)B\eta \rangle \\
&= \langle J\pi (xy^*)B\xi ,B\eta \rangle
= \langle T(xy^*)\xi ,\eta \rangle \\
&= \langle K_T(x,y)\xi ,\eta \rangle ,
\end{split}\end{equation*}
so that $(V,\cK )$ is a Kolmogorov decomposition of $K_T$.
Let $J=J_+-J_-$ be the Jordan decomposition
of $J$ and define
the hermitian kernels
$$K_{\pm }(x,y)=V(x)^{* }J_{\pm }V(y).$$
One can check
that 
$$K_{\pm }(x,\phi (a,y))=K_{\pm }(\phi (a^{-1},x),y)$$
for all $x,y\in \cA $ and $a\in \cU(\cA )$. For $x\in \cA $ define
$$T_{\pm }(x)=K_{\pm }(x,1).$$
Then $T_{\pm }(x)= B^{\*}J\pi (x^*)^{*}J_{\pm }V(1)$
are linear maps on $\cA$ and for $x\in \cA$, $y\in \cU(\cA )$, 
we get
\begin{equation*}\begin{split}
K_{T_{\pm }}(x,y)&=T_{\pm }(xy^{-1})=K_{\pm }(xy^{-1},1) \\
& = K_{\pm }(\phi (y,x),1)=K_{\pm }(x,\phi (y^{-1},1)) \\
& = K_{\pm }(x,y).\end{split}
\end{equation*}
Since $K_{T_{\pm }}$ and $K_{\pm }$ are antilinear in the
second variable and $\cA $ is the linear span of $\cU(\cA )$,
it follows that
$$K_{T_{\pm }}(x,y)=K_{\pm }(x,y)$$
for all $x,y\in \cA$. This implies that
$T_{\pm }$ are disjoint completely positive maps 
such that $T=T_+-T_-$.
The implication (3)$\Rightarrow$(1) follows from Theorem~\ref{baza}.
\end{proof}

Theorem~\ref{witt} suggests how to extend the concept of decomposition 
to arbitrary kernels. In the following we repeatedly use the following
observation: if $L$ is a positive definite
kernel on $X$ and with values in $\cL(\cH)$, and $T\in\cL(\cH)$, 
then the kernel $T^*LT^*$ is positive
definite. Thus, if $K$ is a hermitian kernel, $L$
is a positive definite kernel, both on a
set $X$ and with entries in $\cL(\cH)$, for some Hilbert space $\cH$, 
 then for any $x,y\in X$ 
\begin{equation*}\begin{split}
\left[\begin{array}{cc}
1 & 1 \\
1 & -1 
\end{array}\right]^* &
\left[\begin{array}{cc}
L(x,y) & K(x,y) \\
K(x,y) & L(x,y) 
\end{array}\right] 
\left[\begin{array}{cc}
1 & 1 \\
1 & -1 
\end{array}\right]  \\
 &  =2\left[\begin{array}{cc}
L(x,y)+K(x,y) & 0 \\
0 & L(x,y)-K(x,y) 
\end{array}\right].\end{split}\end{equation*}
Thus, 
$\left[\begin{array}{cc}
L & K \\
K & L
\end{array}
\right]$
is positive definite if and only if both 
$L+K$ and $L-K$ are positive definite, that
is, if and only if the Schwartz condition $-L\leq K\leq L$ holds.

Following U.~Haagerup \cite{Ha}, 
this observation and Theorem~\ref{witt} suggest the 
following definition: a kernel $K:X\times X\rightarrow \lh $ is called
{\it decomposable} if there exist two positive definite
kernels $L_1,L_2:X\times X\rightarrow \lh $
such that the kernel
$\left[\begin{array}{cc}
L_1 & K \\
K^* & L_2
\end{array}
\right]$ is positive definite. Clearly, a kernel $K$ is decomposable if and
only if it is a linear combination of positive definite kernels. Actually, it
is easy to see that any decomposable kernel can be written as linear
combination of at most four positive definite kernels. The next result can be
viewed as an analog of V.~Paulsen dilation theorem \cite{Pa}.

\begin{theorem}\label{decomp}
The kernel $K$ is decomposable if and only if
there is a Hilbert space $\cK$, a mapping 
$V:X\rightarrow \cL(\cH,\cK)$, and a contraction 
$U$ on $\cK$ such that 
$K(x,y)=V(x)^*UV(y)$ for all $x,y\in X$, 
and the set $\{V(x)h\mid x\in X, h\in \cH \}$ is total
in $\cK$.
\end{theorem}

\begin{proof}
If $K(x,y)=V(x)^*UV(y)$ for $x,y\in X$ and some contraction $U$, 
then consider the positive definite kernels
$L_1(x,y)=L_2(x,y)=V(x)^*V(y)$. 
We deduce that

$$
\begin{array}{rcl}
\left[\begin{array}{cc}
L_1(x,y) & K(x,y) \\
K^*(x,y) & L_2(x,y)
\end{array}\right]&=&
\left[\begin{array}{cc}
V(x)^*V(y) & V(x)^*UV(y) \\
V(x)^*U^*V(y) & V(x)^*V(y)
\end{array}\right] \\
& & \\
 &= & \left[\begin{array}{cc}
V(x)^* & 0 \\
0 & V(x)^*
\end{array}\right]
\left[\begin{array}{cc}
I & U \\
U^* & I
\end{array}\right]
\left[\begin{array}{cc}
V(y) & 0 \\
0 & V(y)
\end{array}\right].
\end{array}
$$

\noindent
Since $U$ is a contraction, the matrix $\left[\begin{array}{cc}
I & U \\
U^* & I
\end{array}\right]$
is positive. Next, take $x_1,\ldots ,x_n\in X$ and after 
reshuffling, the matrix
$$\left[\left[\begin{array}{cc}
L_1(x_i,x_j) & K(x_i,x_j) \\
K^*(x_i,x_j) & L_2(x_i,x_j)
\end{array}\right]
\right]_{i,j=1}^n
$$
can be written in the form
$$
\left(\bigoplus _{i=1}^n(V(x_i)^*\oplus V(x_i)^*\right)
\left(\left[\begin{array}{cc}
I & U \\
U^* & I
\end{array}\right]\otimes 
\left[\begin{array}{cccc}
I & I & \ldots & I \\
I & I & \ldots & I \\
\vdots & \vdots & \ddots & \\
I & I & & I 
\end{array}
\right]\right)
\left(\bigoplus _{j=1}^n(V(x_j)\oplus V(x_j)\right),
$$
which shows that the kernel 
$
\left[\begin{array}{cc}
L_1(x,y) & K(x,y) \\
K^*(x,y) & L_2(x,y)
\end{array}\right]
$
is positive definite. 

Conversely, assume that $K$ is decomposable.
We consider the real and imaginary parts of $K$, 
\begin{equation}\label{kaunu}
K_1(x,y)=\frac{1}{2}(K(x,y)+K^*(x,y)),
\end{equation} 
\begin{equation}\label{kadoi}
K_2(x,y)=\frac{1}{2\iac}(K(x,y)-K^*(x,y)),\end{equation}
therefore $K_1$, $K_2$ are hermitian kernels and $K=K_1+\iac K_2$.
Since $K$ is decomposable, there exist positive definite kernels 
$L_1$ and $L_2$ such that the kernel
$
\left[\begin{array}{cc}
L_1 & K \\
K^* & L_2
\end{array}\right]
$
is positive definite.
Therefore, 
\begin{equation*}\begin{split}
\left[\begin{array}{cc}
1 & 1 \\
1 & -1
\end{array}\right]^* &
\left[\begin{array}{cc}
L_1 & K \\
K^* & L_2
\end{array}\right]
\left[\begin{array}{cc}
1 & 1 \\
1 & -1
\end{array}\right] \\
& =\left[\begin{array}{cc}
L_1+K^*+K+L_2 & L_1+K-K^*-L_2 \\ L_1+K^*-K-L_2 & L_1-K^*-K+L_2
\end{array}\right]
\end{split}\end{equation*}
is also a positive definite kernel, which implies that 
$-\frac{1}{2}(L_1+L_2)\leq K_1\leq \frac{1}{2}(L_1+L_2)$.
Similarly, 
\begin{equation*}\begin{split}
\left[\begin{array}{cc}
1 & \iac \\
-\iac  & -1
\end{array}\right]^* &
\left[\begin{array}{cc}
L_1 & K \\
K^* & L_2
\end{array}\right]
\left[\begin{array}{cc}
1 & \iac \\
-\iac & -1
\end{array}\right] \\ &
=\left[\begin{array}{cc} L_1+\iac K^*-\iac K+L_2 & \iac L_1-K^*+K+\iac L_2 \\
-\iac L_1-K+K^*+\iac L_2 & L_1-\iac L_1-\iac K^*+\iac K+L_2 \end{array}\right]
\end{split}\end{equation*}
is a positive definite kernel, which gives that
$-\frac{1}{2}(L_1+L_2)\leq K_2\leq \frac{1}{2}(L_1+L_2)$.
Since $\frac{1}{2}(L_1+L_2)$ is a positive definite kernel, we deduce
from Theorem~\ref{kolmo} that
both $K_1$ and $K_2$ have Kolmogorov decompositions, say
$$K_i(x,y)=V_i(x)^*J_1V_i(y),\quad i=1,2,$$
where $V_i:X\rightarrow \cL(\cH, \cK _i)$ and $J_i$ are fundamental
symmetries on $\cK_i$, $i=1,2$. It follows that 
$$
\begin{array}{rcl}
K(x,y)&=&V_1(x)^*J_1V_1(y)+\iac V_2(x)^*J_2V_2(y) \\
 & & \\
 &=&\left[\begin{array}{c}
V_1(x) \\
V_2(x)
\end{array}
\right]^*
\left[\begin{array}{cc}
J_1 & 0 \\
0 & \iac J_2
\end{array}
\right]
\left[\begin{array}{c}
V_1(y) \\
V_2(y)
\end{array}
\right].
\end{array}
$$
Define 
$
V'(x)=\left[\begin{array}{c}
V_1(x) \\
V_2(x)
\end{array}
\right]
:\cH \rightarrow \cK_1\oplus \cK _2
$ 
and let $\cK$ be the closure in 
$\cK_1\oplus \cK _2$
of the linear space generated by the elements of the 
form $V'(x)h$, $x\in X$ and $h\in \cH$. Finally, define 
$V(x)=P_{\cK }V'(x)$, where $P_{\cK }$ denotes the 
orthogonal projection of $\cK_1\oplus \cK _2$
onto $\cK$. Then the set $\{V(x)h\mid x\in X, h\in \cH \}$
is total in $\cK$ and $K(x,y)=V(x)^*P_{\cK } 
\left[\begin{array}{cc}
J_1 & 0 \\
0 & \iac J_2
\end{array}
\right]
P_{\cK }V(y)$. The operator $U=P_{\cK } 
\left[\begin{array}{cc}
J_1 & 0 \\
0 & \iac J_2
\end{array}
\right]
P_{\cK }$
is a contraction and the proof is concluded.
\end{proof}

\section{Holomorphic kernels}
There are many examples of hermitian kernels which are 
holomorphic on some domain in the complex plane, see for
instance \cite{ADRS}. In all these cases it is known that the 
kernels are associated with reproducing kernel Kre\u\i n (Hilbert) spaces,
and D.~Alpay in \cite{Al} proved a general result in this direction. Our goal
is to extend the result in \cite{Al} to the multi-variable case. In view of
the transcription between Kolmogorov decompositions and reproducing kernel
spaces, e.g.\ see the invariant version in Theorem~\ref{repro}, 
we actually show that the orginal idea of the proof in \cite{Al},
which goes back to \cite{Sc}, can be adapted to this multi-variable setting.

We first review the well-known example of the Szeg\"o kernel
(see \cite{Ar}).
Let $\cG$ be a Hilbert space and denote by $B_r(\xi)$ the open
ball of radius $r$ and center $\xi $,
$B_r(\xi)=\{\eta \in \cG \mid \|\eta -\xi \|<r\}$. We write
$B_r$ instead of $B_r(0)$. For $\xi, \eta \in B_1$, define
\begin{equation}\label{szego}
S(\xi ,\eta )=\frac{1}{1-\langle \eta ,\xi \rangle},\end{equation}
and note that $S$ is a positive definite kernel on 
$B_1$. We now describe its Kolmogorov decomposition. Let
$$F(\cG )=\bigoplus _{n=0}^{\infty }\cG ^{\otimes n},$$
be the Fock space associated to $\cG$, 
where $\cG ^{\otimes 0}=\CC$ 
and $\cG ^{\otimes n}$ is the $n$-fold tensor product
of $\cG$ with itself. Let 
\begin{equation}P_n=(n!)^{-1}\sum _{\pi \in S_n}\hat\pi\end{equation}
be the orthogonal projection of $\cG ^{\otimes n}$ onto its
symmetric part, where 
$$\hat\pi (\xi _1\otimes \ldots \otimes \xi _n)=\xi _{\pi ^{-1}(1)}
\otimes \ldots \otimes \xi _{\pi ^{-1}(n)}$$
and $\pi $ is an element of the permutation group $S_n$
on $n$ symbols.
The symmetric Fock space is 
$F^s(\cG )=(\bigoplus\limits _{n=0}^{\infty }P_n)F(\cG )$.
For $\xi \in B_1$ set $\xi ^{\otimes 0}=1$ and let $\xi ^{\otimes n}$ 
denote the 
$n$-fold tensor product $\xi\otimes \ldots 
\otimes \xi $, $n\geq 1$. Note that
$$\|\bigoplus _{n\geq 0}\xi ^{\otimes n}\|^2=
\sum _{n\geq 0}\|\xi  ^{\otimes n}\|^2=\sum _{n\geq 0}\|\xi \|^{2n}=
\frac{1}{1-\|\xi \|^2}.$$
Hence $\bigoplus\limits _{n\geq 0}\xi ^{\otimes n}\in F^s(\cG )$
and we can define the mapping $V_S$ from $B_1$ into  $F^s(\cG )$, 
\begin{equation}V_S(\xi )=\bigoplus _{n\geq 0}\xi ^{\otimes n},\quad
  \xi\in\cG.\end{equation} 

\begin{lemma}\label{kolmofock}
The pair $(V_S,F^s(\cG ))$ is the Kolmogorov decomposition of the kernel
$S$.
\end{lemma}

\begin{proof}
$V_S(\xi )$ is also viewed as a bounded linear operator from 
$\CC $ into $F^s(\cG )$ by $V_S(\xi )\lambda =\lambda V_S(\xi )$,
$\lambda \in \CC $, so that, for $\xi ,\eta \in B_1$, 

\begin{equation*}\begin{split}
V_S(\xi )^*V_S(\eta )&=\langle V_S(\eta ),V_S(\xi )\rangle \\
 &=\sum _{n\geq 0}\langle \eta ^{\otimes n},\xi ^{\otimes n}\rangle \\
 &=\sum _{n\geq 0}\langle \eta ,\xi \rangle ^n=
\displaystyle\frac{1}{1-\langle \eta ,\xi \rangle}=S(\xi ,\eta ).
\end{split}\end{equation*}
The set $\{V_S(\xi )\mid \xi \in B_1\}$ is total in 
$F^s(\cG )$ since for $n\geq 1$ and $\xi\in\cG$ we have
$\displaystyle\frac{d^n}{dt^n}V(t\xi )|_{t=0}=n!\xi ^{\otimes n}.$
\end{proof}

The reproducing kernel Hilbert space
associated to the Szeg\"o kernel $S$, see
\eqref{szego}, is given by the completion of the linear space generated 
by the functions $s_{\eta }=S(\cdot ,\eta )$, $\eta \in B_1$, 
with respect to the inner product defined by
$$\langle s_{\eta },s_{\xi }\rangle=S(\xi ,\eta ).$$
Note that there exists a unitary operator $\cF $ from 
the reproducing kernel Hilbert space associated to the Szeg\"o kernel $S$ 
onto $F^s(\cG )$ such that
$\cF s_{\xi }=V_S(\xi )$, $\xi \in B_1$.

We now explore the fact that $S$ is 
a holomorphic kernel. We use the terminology and results from
\cite{Mu}, \cite{Di} for holomorphic functions in infinite dimensions.
Thus, we say that a function $f$ defined on the open subset $\cO$ 
of $\cG$ is holomorphic if $f$ is continuous on $\cO$ and 
for all $\eta \in \cO$, $\xi \in \cG$,
the mapping $\lambda \rightarrow f(\eta +\lambda \xi )$
is holomorphic on the open set $\{\lambda \in \CC \mid 
\eta +\lambda \xi \in \cO\}$. One easily check that $V_S$
is holomorphic on $B_1$, 
therefore $S(\xi ,\cdot )$ is holomorphic on $B_1$ for each
fixed $\xi \in B_1$. 
We also notice that the reproducing kernel Hilbert space
associated to $S$ consists of anti-holomorphic functions on
$B_1$. It is somewhat more convenient to replace this space
by a Hilbert space of holomorphic functions on $B_1$.
Thus, define the holomorphic function $a_{\xi }=S(\xi ,\cdot )$ on $B_1$
for each $\xi \in B_1$,
and let $H^2(\cG )$ denote the completion of 
the linear space generated 
by the functions $a_{\xi }$, $\xi \in B_1$, 
with respect to the inner product defined by
$$\langle a_{\xi },a_{\eta }\rangle=S(\xi ,\eta ),\quad \xi,\eta\in B_1.$$
We notice that $H^2(\cG )$ is an anti-unitary copy of the 
reproducing kernel Hilbert space of $S$.

We say that a hermitian kernel $K$, defined on an open subset $\cO$ of $\cG$,
is {\it holomorphic on} $\cO$ if $K(\xi ,\cdot )$ is holomorphic on $\cO$ for
each fixed $\xi \in \cO$.

\begin{theorem}\label{alpay}
Let $K$ be a holomorphic hermitian kernel on $B_r$, with $r>0$, and valued in
$\cL(\cG)$ for some Hilbert space $\cG$. Then there exists $0<r'\leq r$ and
a Kolmogorov decomposition $(V;\cK)$ of $K|B_{r'}\times B_{r'}$, such that
$V\colon B_{r'}\ra\cL(\cG,\cK)$ is holomorphic.
\end{theorem}

\begin{proof}
Let $K$ be a holomorphic kernel on $B_r$, $r>0$. Since $K$
is hermitian, it follows that $K(\cdot ,\eta )$ is anti-holomorphic
on $B_r$ for each $\eta \in B_r$. It is convenient to reformulate
this fact as follows. Let $\{e_{\alpha }\}_{\alpha \in A}$
be an orthonormal basis for $\cG$. We define the mapping 
$$\xi =\sum _{\alpha \in A}\langle \xi ,e_{\alpha }\rangle e_{\alpha }
\rightarrow 
\sum _{\alpha \in A}\overline{\langle \xi ,e_{\alpha }\rangle }e_{\alpha }=
\xi ^*,
$$
so that the function $f(\xi ,\eta )=K(\xi ^*,\eta )$
is separately holomorphic on $B_r\times B_r$. By Hartogs' Theorem
(\cite{Mu}, Theorem~36.8), $f$ is holomorphic on $B_r\times B_r$.
By (\cite{Mu}, Proposition~8.6), $f$ is locally bounded.
We first suppose that $r>1$. Hence
there exist $1<\rho <r$ and $C>0$ such that:
\begin{equation}\label{cinciunu}
|K(\xi ,\eta )|\leq C \quad \mbox{for}\quad \xi ,\eta \in B_{\rho }
\end{equation}
and 
\begin{equation}\label{cincidoi}
K(\xi ^*,\eta )=\sum _{m\geq 0}p_m(\xi, \eta)
\end{equation}
uniformly on $B_{\rho }$, where each $p_m$, $m\geq 0$, is an 
$m$-homogeneous continuous polynomial on $\cG \times \cG $.
That is (see \cite{Mu} or \cite{Di}, Chapter ~1), there exists 
a continuous linear mapping $A_m$ on $P_m((\cG \times \cG)^{\otimes m})$
such that 
\begin{equation}\label{cincitrei}
p_m(\xi ,\eta )=A_m((\xi ,\eta )^{\otimes m})
\end{equation}
for all $\xi ,\eta \in \cG $. 

Using Cauchy Inequalities, \cite{Di}, Proposition~3.2, 
for $B_{\rho }$, we deduce
\begin{equation}\label{cincipatru}
\|A_m\|\leq \|p_m\|\leq C\displaystyle\frac{1}{\rho ^m},
\end{equation}
hence
\begin{equation}\label{cincicinci}
\sum _{m\geq 0}\|A_m\|^2\leq C\sum _{m\geq 0}
\displaystyle\frac{1}{\rho ^{2m}}=C
\displaystyle\frac{1}{1-1/\rho ^2}=C'<\infty.
\end{equation}

The previous facts are valid with respect to the norm
$\|(\xi ,\eta )\|=\max \{\|\xi \|,\|\eta \|\}$. 
Since this norm is equivalent to the Hilbert norm 
$\|(\xi ,\eta )\|=\sqrt{\|\xi \|^2+\|\eta \|^2}$, 
we deduce that each $A_m$ is also continuous with respect
to this Hilbert norm 
on $\cG \times \cG $. By Riesz representation theorem, there exist
$a_m\in P_m(\cG \times \cG )^{\otimes m}$, $m\geq 0$, such that 
\begin{equation}\label{cincisase}
A_m((\xi ,\eta )^{\otimes m})=\langle (\xi ,\eta )^{\otimes m},a_m\rangle
\end{equation}
and 
\begin{equation}\label{cincisapte}
\|a_m\|=\|A_m\|
\end{equation}
(with $a_0=A_0\in \CC$).
Since $P_m(\cG \times \cG )^{\otimes m}$ is isometrically isomorphic 
to $(P_m\cG ^{\otimes m})^{\oplus (m+1)}$, we deduce that 
there are $a_m^k\in P_m\cG ^{\otimes m}$, $k=0,\ldots ,m$, such that
\begin{equation}\label{cinciopt}
\langle (\xi ,\eta )^{\otimes m},a_m\rangle=
\sum _{k=0}^m\langle b_m^k(\xi ,\eta ),a_m^k\rangle
\end{equation}
and 
\begin{equation}\label{cincinoua} 
\sum _{k=0}^m\|a_m^k\|^2=\|a_m\|^2,
\end{equation}
where $b_0^0=1$ and $b_m^k(\xi,\eta)=\xi ^{\otimes (m-k)}\otimes
\eta ^{\otimes k}$, $m\geq 1$, $k=0,\ldots ,m$.

We now show that for each fixed $\xi \in B_1$, $g_{\xi }(\eta )=
K(\xi ,\eta )$ belongs to $H^2(\cG )$.
By \eqref{cincidoi}, 
\eqref{cincitrei}, \eqref{cincisase}, and \eqref{cinciopt},
$$K(\xi ,\eta )=\sum _{m\geq 0}\sum _{k=0}^m
\langle b_m^k(\xi ^*,\eta ),a_m^k\rangle ,$$
and the series converges absolutely on $\eta $ by \eqref{cincipatru}.
Reordering to $m$-homogeneous terms in $\eta $, 
$$K(\xi ,\eta )=\sum _{k\geq 0}\sum _{m\geq k}
\langle b_m^k(\xi ^*,\eta ),a_m^k\rangle .$$
For fixed $\xi$ define $F_k(\eta )=\sum _{m\geq k}
\langle b_m^k(\xi ^*,\eta ),a_m^k\rangle $. By Schwarz inequality,
\begin{equation*}\begin{split}
\|F_k(\eta )\|^2&\leq (\sum _{m\geq k}\|b_m^k(\xi ^*,\eta )\|^2)
(\sum _{m\geq k}\|a_m^k\|^2) \\
 &=\|\eta \|^{2k}
(\sum _{m\geq k}\|\xi \|^{2(m-k)})
(\sum _{m\geq k}\|a_m^k\|^2) \\
 &=\displaystyle\frac{\|\eta \|^{2k}}{1-\|\xi \|^2}
\sum _{m\geq k}\|a_m^k\|^2,
\end{split}\end{equation*}
which implies 
$$\|F_k\|\leq 
\displaystyle\frac{1}{1-\|\xi \|^2}
\sum _{m\geq k}\|a_m^k\|^2.$$ 
Finally, 
$$\|g_{\xi }\|_{H^2(\cG )}^2=\sum _{k\geq 0}\|F_k\|^2
\leq \displaystyle\frac{1}{1-\|\xi \|^2}
\sum _{k\geq 0}\sum _{m\geq k}\|a_m^k\|^2.$$
By \eqref{cincicinci},
\eqref{cincisapte}, and 
\eqref{cincinoua},
we deduce that 
$$\|g_{\xi }\|_{H^2(\cG )}^2\leq 
\displaystyle\frac{1}{1-\|\xi \|^2}
\sum _{m\geq 0}\|A_m\|^2\leq C'\displaystyle\frac{1}{1-\|\xi \|^2}.$$
This shows that $g_{\xi }\in H^2(\cG )$ and, more than that, the 
formula 
$$\cP a_\xi =g_{\xi }$$
gives a bounded linear operator $\cP $ on $H^2(\cG )$, such that
$$K(\xi ,\eta )=g_{\xi }(\eta )=(\cP a_{\xi })(\eta )=
\langle \cP a_{\xi },a_{\eta }\rangle _{H^2(\cG )}.$$
This implies that $\cP $ is selfadjoint and let $\cP=\cP_+-\cP_-$
be its Jordan decomposition, where $\cP _{\pm }$ are positive operators
on $H^2(\cG )$. Then $K$ can be written as the difference 
of two positive definite kernels. By \cite{Sc} and Theorem~2.3, 
$K$ has a Kolmogorov decomposition, and it is easy to see that this Kolmogorov
decomposition has the required holomorphy property. In case $r\leq 1$, 
a scalling argument as in \cite{Al} concludes the proof.
\end{proof}

As mentioned, the above proof is based on the same idea as in \cite{Al}, 
which goes back to \cite{Sc}. An interesting aspect of this idea is 
that once again the Szeg\"o kernel
$S$ has a certain universality property, that is, any 
holomorphic kernel on $B_r$, $r>1$, 
is the image of $S$ through a bounded selfadjoint operator on 
$H^2(\cG )$. A different kind of universality property 
of $S$, related to the solution of the Nevanlinna-Pick 
interpolation problem, was established in \cite{AC}.

Finally we apply Theorem~\ref{alpay} to show that non-hermitian holomorphic
kernels are decomposable. A kernel $K\colon \cO\times\cO\ra\cL(\cG)$, where
$\cO$ is an open subset of some Hilbert space $\cG$, 
is {\em holomorphic} on $\cO$ 
if, for any fixed $\xi\in\cO$, the function $K(\xi,\cdot)\colon\cO\ra\cL(\cG)$
is holomorphic and,
for any fixed $\eta\in\cO$, the function $K(\cdot,\eta)\colon\cO\ra\cL(\cG)$ 
is anti-holomorphic.

\begin{corollary} Let $K$ be a holomorphic kernel on $B_r$, 
with $r>0$, and valued in
$\cL(\cG)$ for some Hilbert space $\cG$. Then there exists $0<r'\leq r$, a
Hilbert space $\cK$, a holomorphic mapping 
$V:B_{r'}\rightarrow \cL(\cG,\cK)$, and a contraction $U$ on $\cK$, such that 
$K(\xi,\eta)=V(\xi)^*UV(\eta)$ for all $\xi,\eta\in B_{r'}$,
and the set $\{V(\xi)h\mid \xi\in B_{r'}, h\in \cG \}$ is total in $\cK$.
\end{corollary}

\begin{proof} We consider the real part $K_1$ \eqref{kaunu} and, respectively,
the imaginary part $K_2$ \eqref{kadoi} of $K$ and note that both are
holomorphic hermitian kernels. Then we apply Theorem~\ref{alpay} to produce
holomorphic Kolmogorov decompositions of $K_1$ and $K_2$ on a possibly
smaller, but nontrivial, ball $B_{r'}$ in $\cG$ and, proceeding as in
the proof of Theorem~\ref{decomp}, we obtain the decomposition of $K$ as
required.\end{proof}

\end{document}